\documentclass[a4paper, 12pt]{article}
\usepackage[top=2cm, bottom=2.5cm, left=2.2cm, right=2.2cm]{geometry}
\usepackage{amsmath,amsthm,amssymb}

\usepackage{color}
\usepackage{makeidx}
\usepackage{hyperref}
\makeindex

\numberwithin{equation}{section}

{\theoremstyle{definition}\newtheorem{definition}{Definition}[section]
\newtheorem{defnot}[definition]{Definitions and Notation}

\newtheorem{remark}[definition]{Remark}
\newtheorem{remarks}[definition]{Remarks}
\newtheorem{facts}[definition]{Facts}

}
\newtheorem{proposition}[definition]{Proposition}
\newtheorem{proposition-definition}[definition]{Proposition-Definition}
\newtheorem{lemma}[definition]{Lemma}
\newtheorem{theorem}[definition]{Theorem}
\newtheorem{corollary}[definition]{Corollary}

\newcommand{\R}{\mathbb{R}}
\newcommand{\N}{\mathbb{N}}

\newcommand{\cF}{\mathcal{F}}
\newcommand{\cG}{\mathcal{G}}

\newcommand{\cJ}{\mathcal{J}}
\newcommand{\cE}{\mathcal{E}}

\newcommand{\cL}{\mathcal{L}}

\newcommand{\cN}{\mathcal{N}}
\newcommand{\cK}{\mathcal{K}}
\newcommand{\C}{\mathbb{C}}

\newcommand{\cA}{\mathcal{A}}

\newcommand{\Z}{\mathbb{Z}}
\newcommand{\cP}{\mathcal{P}}

\newcommand{\cH}{\mathcal{H}}
\newcommand{\cM}{\mathcal{M}}

\newcommand{\cS}{\mathcal{S}}

\newcommand{\ie}{{i.e.}\/ }
\newcommand{\eg}{{\it e.g.}\/ }
\newcommand{\cf}{{\it cf.}\/ }

\newcommand{\gA}{\mathfrak{A}}
\newcommand{\gJ}{\mathfrak{J}}

\def\gpd{\,\lower1pt\hbox{$\longrightarrow$}\hskip-.24in\raise2pt
             \hbox{$\longrightarrow$}\,}
             
\newcommand{\BIG}{gauge adiabatic groupoid}
\newcommand{\BIF}{gauge adiabatic foliation}
\newcommand{\boinG}{G_{ga}}

\everymath{\displaystyle}

\reversemarginpar

\usepackage{color}

\begin{document}

\renewcommand\theenumi{\alph{enumi}}
\renewcommand\labelenumi{\rm {\theenumi})}
\renewcommand\theenumii{\roman{enumii}}

\begin{center}
{\Large\bf Adiabatic groupoid, crossed product by $\R_+^*$ and Pseudodifferential calculus

\bigskip

{\sc by Claire Debord and Georges Skandalis}
}
\end{center}

{\footnotesize
Laboratoire de Math\'ematiques, UMR 6620 - CNRS
\vskip-4pt Universit\'e Blaise Pascal, BP {\bf 8002}
\vskip-4pt F-63171 Aubi\`ere cedex, France
\vskip-4pt claire.debord@math.univ-bpclermont.fr

\vskip 2pt Institut de Math{\'e}matiques de Jussieu, (UMR 7586), Universit\'e Paris Diderot (Paris 7)
\vskip-4pt  UFR de Math\'ematiques, {\sc CP} {\bf 7012} - B\^atiment Sophie Germain 
\vskip-4pt  5 rue Thomas Mann, 75205 Paris CEDEX 13, France
\vskip-4pt skandalis@math.univ-paris-diderot.fr
}
\bigskip

\centerline{\bf Abstract}

We consider the crossed product $G_{ga}$ of the natural action of $\R_+^*$ on the adiabatic groupoid $G_{ad}$ associated with any Lie groupoid $G$. We construct an explicit Morita equivalence between the exact sequence of order $0$ pseudodifferential operators on $G$ and (a restriction of) the natural exact sequence associated with $G_{ga}$. As an important intermediate step, we express a pseudodifferential operator on $G$ as an integral associated to a smoothing operator on the adiabatic groupoid $G_{ad}$ of $G$.
\tableofcontents
\setcounter{section}{-1}

\section{Introduction}

Smooth groupoids are intimately linked to pseudodifferential calculi. Indeed, to every smooth groupoid is naturally associated a pseudodifferential calculus and therefore an analytic index. Furthermore many pseudodifferential calculi have been shown to be the ones associated to naturally defined groupoids. The groupoid approach may then give a natural geometric insight to these calculi and the corresponding index theorems.  See \cite{DebordLescure} for an overview on the subject.

The use of groupoids in relation with index theory may be traced back to \cite{ConnesLNM} where A. Connes introduced the longitudinal pseudodifferential calculus on a foliation $(M,F)$. He thus constructs an extension $\Psi^*(M,F)$ of the foliation $C^*$-algebra $C^*(M,F)$ which gives rise to an exact sequence $$0\to C^*(M,F)\to \Psi^*(M,F)\to C(S^*F)\to 0$$  It appeared quite naturally that Connes' construction only used the (longitudinal) smooth structure of the holonomy groupoid and could therefore be extended to any smooth (longitudinally) groupoid (\cf \cite{MonthPie, NWX, LMN}).

Some previously defined pseudodifferential calculi were recognized as being the ones associated with natural groupoids as for instance the groupoid defined by B. Monthubert in \cite{MonthT1, MonthT2} was shown to be suitable for R.B. Melrose's $b$-calculus \cite{MelroseLivre}. Moreover, it appeared in a work of J-M. Lescure (\cite{Lescure}) that this calculus is the natural calculus associated to conical pseudo-manifolds.

In \cite{ConnesNCG}, A. Connes showed that the analytic index on a compact manifold can in fact be described in a way not involving (pseudo)differential operators at all, just by using a construction of a deformation groupoid, called the ``tangent groupoid''. This idea was used in  \cite{HilsumSk}, and extended in \cite{MonthPie} to the general case of a smooth groupoid, where the authors associated to every smooth groupoid $G$ an \emph{adiabatic groupoid}, which is obtained applying the ``deformation to the normal cone'' construction to the inclusion $G^{(0)}\to G$ of the unit space of $G$ into $G$. The groupoid constructed in this way is the union $G_{ad}=G\times \R^*\cup \gA G\times \{0\}$ endowed with a natural smooth structure, where $\gA G$ is the total space of the algebroid of $G$ \ie of the normal bundle to the inclusion $G^{(0)}\to G$.  They then showed that the connecting map of the corresponding exact sequence $$0\to C_0(\R_
+^*)\otimes C^*(G)\to C^*(G_{ad}^+)\to C^*(\gA G)\simeq C_0(\gA^* G)\to 0\eqno (1)$$ is the analytic index, where $G_{ad}^+=G\times \R^*_+\cup \gA G\times \{0\}$ is the restriction of $G_{ad}$ over $G^{(0)}\times \R_+$.

In the present paper, extending ideas of Aastrup, Melo, Monthubert and Schrohe \cite{AMMS}, we go one step further in this direction, showing that the (order $0$) pseudodifferential operators on a smooth groupoid can also be described as convolution operators by smooth functions on a suitable groupoid. The groupoid that we use is the crossed product of the adiabatic groupoid by the natural action of the group $\R_+^*$. Since exact sequence (1) is equivariant with respect to this action, we find an exact sequence $$0\to C^*(G)\otimes \cK\simeq\Big(C_0(\R_+^*)\otimes C^*(G)\Big)\rtimes \R_+^*\to C^*(G_{ad}^+)\rtimes \R_+^*\to C_0(\gA^* G)\rtimes \R_+^*\to 0.$$
The algebra $C_0(\gA^* G)\rtimes \R_+^*$ naturally contains $C(S^*\gA G)\otimes \cK$ as an ideal. The main result of this paper is that the corresponding exact sequence $$0\to C^*(G)\otimes \cK\to J(G)\rtimes \R_+^*\to C(S^*\gA G)\otimes \cK\to 0$$ is related to the exact sequence of pseudodifferential operators $$0\to C^*(G)\to \Psi^*(G)\to C(S^*\gA G)\to 0$$ via a Morita equivalence.

Our groupoid is the generalization of the one constructed in \cite{AMMS} for the case of the ordinary pseudodifferential calculus on a compact manifold $M$. It was shown there that the algebra associated to this groupoid is isomorphic to the algebra of Green operators (of order $0$ and class $0$) in the Boutet de Monvel calculus (\cf \cite{Grubb, Schrohe}) which in turn is known to be Morita equivalent to $\Psi^*(M)$. The proof in \cite{AMMS} is somewhat indirect, using Voiculescu's theorem to prove that two exact sequences coincide. 

Our proof here is much more direct and therefore extends immediately to the general groupoid case: we explicitly construct a bimodule $\cE$ which is a Morita equivalence between the algebras $\Psi^*(G)$ and $J(G)\rtimes \R_+^*$. As an important intermediate step, we express a pseudodifferential operator on $G$ as an integral associated to a smoothing operator on the adiabatic groupoid $G_{ad}$ of $G$ (Theorem \ref{pseudosintegrales}). We should point out that J-M Lescure had previously observed that pseudodifferential operators on $\R^n$ arise as integrals of some functions on the tangent groupoid of $\R^n$ (private communication).

Next, we show that the $C^*$-algebra $J(G)\rtimes \R_+^*$ is stable and therefore isomorphic to $\Psi^*(G)\otimes \cK$, and the Hilbert module $\cE$ is isomorphic to Kasparov's absorbing module $\cH_{\Psi^*(G)}$. To that end, we use the fact that positive order pseudodifferential operators define regular operators which was established by S. Vassout in \cite{chief}. Actually, using pseudodifferential calculus of complex order and complex powers also from \cite{chief}, we construct an isomorphism of $J(G)$ onto a crossed product $\Psi^*(G)\rtimes \R$ which intertwines the action $\alpha $ of $\R_+^*$ on $J(G)$ with the dual action on the crossed product $\Psi^*(G)\rtimes \R$.

In a forthcoming paper (\cite{DSBM}), we will actually show that the algebra $J(G)\rtimes \R_+^*$ is not only isomorphic, but in fact \emph{equal} to the algebra of Green operators - generalized to any Lie groupoid.

The paper is organized as follows:
\begin{itemize}
\item In the first section, we recall the construction of the deformation to the normal cone and the space of Schwartz functions on it. We also characterize the space of those functions whose Fourier transform vanishes at infinite order at $0$.
\item In the second section we recall some definitions concerning Lie groupoids, and in particular the construction of the associated ``adiabatic'' groupoid and its crossed product by the canonical action of $\R_+^*$.
\item Finally, in the last two sections, we construct the  bimodule which is used for the Morita equivalence: in the third section we construct the smooth module $\cE^\infty$, whose $C^*$-completion the Hilbert $C^*$-module $\cE$ is shown to indeed define the desired Morita equivalence in section 4.
\end{itemize}

\textbf{Remark on the general setting.} \emph{For the simplicity of our exposition, we will assume that the groupoids involved are smooth and Hausdorff. We will further in some places reduce to the case where the unit space is compact. Of course one easily extends our constructions and results to more general settings: longitudinally smooth groupoids on manifolds with corners, continuous family groupoids (\cf \cite{Paterson}), non Hausdorff case (\cf \cite{ConnesSurvey}), and even groupoids associated with singular foliations (\cf \cite{AndrSk1})... To proceed, one only has to consider the appropriate spaces of functions on the groupoids - that appear in the above cited papers.\\
Furthermore, in what follows, $C^*(G)$ may be either the reduced or full $C^*$-algebra of the groupoid $G$. The choice is left to the reader!}

\medskip \textbf{Acknowledgements.} The first author would like to thank Jean-Marie Lescure for several illuminating discussions and his valuable help for the understanding of pseudodifferantial operators. She would also like to thank the CNRS and the IMJ for hosting her during the spring term and especially the chief of the operator algebra team for the wonderful work conditions offered.

\section{Schwartz spaces on deformations to normal cones}

The adiabatic groupoid which is the main ingredient in our construction is a special case of a geometric object called \emph{deformation to the normal cone.}

\subsection{Deformation to a normal cone}\label{DNC}

In this section, we recall this construction and define some function spaces on it: the associated space of functions of Schwartz decay \emph{alla} P. Carrillo-Rouse (\cite{CR}) and some subspaces which will be needed in the sequel of the paper.

Let $M_0$ be a smooth compact submanifold of a smooth manifold $M$ with normal bundle $\cN$. As a set, the deformation to the normal cone is $D(M_0,M)=M\times \R^*\cup \cN\times \{0\}$. 
In order to recall its smooth structure, we fix an exponential map, which is a diffeomorphism $\theta$ from a neighborhood $V'$ of the zero section $M_0$ in $\cN$ to a neighborhood $V$ of $M_0$ in $M$. 

We may cover $D(M_0,M)$ with two open sets $V\times \R^*$ and $W=\cN\times \{0\}\cup V\times \R^*$: we endow $D(M_0,M)$ with the smooth structure for which the map $\Theta:(x,X,t)\mapsto ( \theta(x,tX),t)$ (for $t\ne 0$) and $\Theta:(x,X,0)\mapsto (x,X,0)$ is a diffeomorphism from  $W'=\{(x,X,t)\in \cN\times \R;\ (x,tX)\in V'\}$ to $W$.

The group $\R_+^*$ acts smoothly on $D(M_0,M)$: for $t\in \R_+^*$ put $\alpha_t(z,\lambda)=(z,t\lambda)$ for $z\in M$ and $\lambda \in \R^*$ and $\alpha _t(x,U,0)=(x,\frac{U}t,0)$ for $x\in M_0$ and $U\in \cN_x$.

\subsection{Schwartz decay on $\R_+^*\times V$} \label{renee}

If $V$ is a smooth non necessarily compact manifold, define $\cS(\R_+^*;C_c^\infty(V))$ to be the space of smooth functions $t\mapsto f_t$ from $\R_+^*$ to $C_c^\infty(V)$ such that all the $f_t$ have support on a given compact subset of $V$ and $t\mapsto f_t$ has rapid decay with respect to all natural norms of $C_c^\infty(V)$. Equivalently $(f_t)$ is such that  the function $g:V\times \R\to \C$ defined by $g(z,t)=f_{\frac{t}{1-t}}(z)$ for $0<t<1$ and $g(z,t)=0$ otherwise is smooth with compact support. 

We will also consider $\cS(\R^*;C_c^\infty(V))$ which is the space of smooth functions $t\mapsto f_t$ from $\R^*$ to $C_c^\infty(V)$ such that all the $f_t$ have support on a given compact subset of $V$ and $t\mapsto f_t$ has rapid decay with respect of all norms at $0$ and $\pm\infty$.

\subsection{Schwartz functions on a vector bundle}

\label{Schwartzfibre}There is a description of the space of Schwartz functions on the total space of a vector bundle in \cite{CR}. Here are two equivalent ways to defining this algebra in a way that it is obvious it is independent of choices of charts.

\begin{definition}
Let $E$ be a smooth real vector bundle over a smooth compact manifold $M$.\begin{enumerate}
\item  Consider $E$ as an open subspace of the bundle of spheres $S_E$ where $M_\infty$ is the set of points at infinity. The algebra $\cS(E)$ is the space of smooth functions on $S_E$ vanishing at any point of $M_\infty$ as well as all its derivatives.

\item For $t\in \R_+^*$, let $\beta_t$ denote the map $(x,\xi)\mapsto (x,t\xi)$ from $E$ to $E$. Then define $\cS_\beta(E)$ to be the space of functions of the form $z\mapsto \int_0^{+\infty}f_t(\beta_t(z))\,dt$ where $f\in \cS(\R_+^*;C_c^\infty(E))$.
\end{enumerate}
\end{definition}

The fact that $\cS_\beta(E)=\cS(E)$ is quite obvious. Indeed if $f\in \cS(\R_+^*;C_c^\infty(E))$, then one defines $g\in C_c^\infty(S_E\times \R)$ by $g(z,t)=f_{\frac{t}{1-t}}(\beta_{\frac{t}{1-t}}(z))$ if $0<t<1$ and $z\in E$ and  $g(z,t)=0$ elsewhere. Then one may integrate it and obtain that $z\mapsto \int_0^{+\infty}f_t(\beta_t(z))\,dt$ is smooth on $S_E$ and vanishes as well as all its derivatives on $M_\infty$.

Conversely, we have to show that the map $\varphi _\beta:\cS(\R_+^*;C_c^\infty(E))\to \cS(E)$ is onto, where $\varphi _\beta(f):z\mapsto \int_0^{+\infty}f_t(\beta_t(z))\,dt$. Choose a metric on $E$; if $g\in \cS(E)$, we may put $f(z,t)=\frac{h(\|z\|^2+t^2)}tg(\beta_{t^{-1}}(z))$ where $h\in C_c^{\infty}(\R)$ vanishes near $0$ and $\int_0^{+\infty}h(s^2)\frac {ds}s=1$. Then $\varphi_\beta(f)(z)=\int_0^{+\infty}h(t^2(1+\|z\|^2))g(z)\,\frac{dt}t=g(z)$.

\begin{remark}
The Fourier transform is an isomorphism of $\cS_\beta(E)$ with $\cS_{\beta^*}(E^*)$. 
\end{remark}

We will use the following rather easy result:

\begin{proposition}\label{J(G)1} Let $f=(f_t)_{t\in \R}\in \cS (E\times \R)$\begin{enumerate}
\item For all $g\in C_c^{\infty}(E)$ the function $F:M\times \R\to \C$ defined by $F(x,0)=\hat f_0(x,0) g(x,0)$ and $F(x,t)=t^{-p}\int_{E_x} g(x,U)f_t\circ \theta (x,U)\,dU$ for $t\ne 0$  is smooth on $M\times \R$.
\item The following are equivalent:\begin{enumerate}\renewcommand\theenumii{\roman{enumii}}
\renewcommand\labelenumii{\rm ({\theenumii})}
\item For all $g\in \cS(E)$ the function $t\mapsto \int_E g(x,U)f (x,\frac Ut,t)\,dx\,dU$ vanishes as well as all its derivatives at $0$.
\item For all $g\in \cS(E)$ the (smooth) function $(x,t)\mapsto \int_{E_x} g(x,U)f (x,\frac Ut,t)\,dU$ vanishes as well as all its derivatives on $M\times \{0\}\subset M\times \R$.
\item The function $(x,\xi,t)\mapsto \hat f_t(x,\xi)$ vanishes as well as all its derivatives on $M\times \{0\}$ sitting in $E^*\times \R$ as zero section.
\end{enumerate}
\end{enumerate}

\begin{proof} Parseval's formula yields \ $t^{-p}\int_{E_x} g(x,U)f (x,\frac Ut,t)\,dU=c\int _{E_x^*}\hat g(x,-\xi)\hat f_t(x,t\xi)\,d\xi$ (where $c$ is a suitable constant and $p$ is the dimension of $E$). 

\begin{enumerate}
\item The function $(x,\xi,t)\mapsto \hat g(x,-\xi)\hat f_t(x,t\xi)$ lies in $\cS(E^*\times \R)$, thus (a) follows.

\item (ii)$\Rightarrow$(i) is obvious. Conversely if (i) is satisfied, writing $\int_{E_x} g(x,U)f (x,\frac Ut,t)\,dU=t^kh_k(x)+o(t^k)$, and applying (i) to $g_1(x,U)=g(x,U)\overline{h_k(x)}$, we find $h_k=0$; thus by induction, we get (ii).

We may write the Taylor expansion $\hat f_t(x,t\xi)=\sum_{j=0}^kt^ja_j(x,\xi)+t^{k+1}R(x,\xi,t)$ where $a_j$ are polynomials in $\xi$ (of degree $\le j$). It follows that $\int_{E_x} g(x,U)f (x,\frac Ut,t)\,dU=o(t^{k+p})$ if and only if $a_j(x,\xi)=0$, whence (ii)$\iff$(iii).
\qedhere
\end{enumerate}
\end{proof}
\end{proposition}

\begin{remark}
It is natural in this proposition to consider $f_t(x,U)$ as densities on $E_x$ rather than functions and therefore introduce a factor $t^{-p}$ in the integrals (i) and (ii). This factor of course has no incidence on the proposition...
\end{remark}

\subsection{Trivial deformation to the normal cone}

\label{actionAlpha} We consider here $D(M,E)$ where the manifold $M$ is sitting as zero section in the total space $E$ of a vector bundle over $M$. Then $D(M,E)=E\times \R$.
Here, $\R_+^*$ acts on $E\times \R$ by $\alpha_t(x,U,\lambda)=(x,\frac{U}{t},t\lambda)$. 

For $f\in \cS(\R_+^*;C_c^{\infty}(E\times \R))$, put $\varphi_\alpha(f)(z)=\int_0^{+\infty}f_t(\alpha_t(z))\,dt$. The image $\cS_\alpha (E\times \R)$ of $\varphi_\alpha$ is the set of $g\in \cS(E\times \R)$ such that the function $(x,U,t)\mapsto \|tU\|$ is bounded on the support of $g$.

Indeed, if $f\in \cS(\R_+^*;C_c^{\infty}(E\times \R))$, one checks immediately the support requirement and it is quite easy to check locally that $\varphi_\alpha(f)\in \cS(E\times \R)$. Conversely, let $g\in \cS(E\times \R)$ with the support requirements; take $\chi\in C_c^\infty(\R)$ equal to $1$ near $0$ and $h\in C_c^{\infty}(\R)$ a function which vanishes near $0$ and satisfies $\int_0^{+\infty}h(s^2)\frac {ds}s=1$ as previously. We may set $f_1(x,U,\lambda,t)=\frac{h(\lambda ^2+t^{2})}tg(x,tU,\frac{\lambda}{t})(1-\chi(\frac{\lambda^2}{t^2}))$ and $f_2(x,U,\lambda,t)=\frac{h(\|U\|^2+t^{-2})}tg(x,tU,\frac{\lambda}{t})\chi(\frac{\lambda^2}{t^2})$.\\
Note that for any $a>0$, the function $t\mapsto  \frac{1}tg(x,tU,\frac{\lambda}t)$ obviously belongs to $\cS(\R_+^*;C_c^{\infty}(E_{\geq a}\times \R_{\geq a}))$ where $E_{\geq a}=\{U\in E \ ; \|U\|\geq a\}$ and $\R_{\geq a}=\R\setminus ]-a,a[$.

For small enough $|\lambda|$ either $h(\lambda ^2+t^{2})=0$ or $\chi(\frac{\lambda^2}{t^2})=1$ thus $f_1$ vanishes. For big enough $|\lambda|$ or $t$, $h(\lambda ^2+t^{2})=0$. Moreover $f_1$ has rapid decay when $t\rightarrow 0$. Finally, $f_1$ belongs to $\cS(\R_+^*;C_c^{\infty}(E\times \R))$.
Similarly $f_2$ vanishes for $t$ near $ 0 $ or for small enough $\|U\|$ and $t\rightarrow \infty$. Anyway $f_2$ has rapid decay when $t\rightarrow \infty$ and thus $f_2$ belongs to $\cS(\R_+^*;C_c^{\infty}(E\times \R))$. One can easily check that $\varphi _\alpha (f_1)+\varphi _\alpha (f_2)=g$.

\subsection{Schwartz functions on a deformation to the normal cone}

In the same way as for bundles, we define $\cS_\alpha(D(M_0,M))$ to be the set of integrals $\varphi _\alpha(f):z\mapsto \int_0^{+\infty}f_t(\alpha_t(z))\,dt$ where $t\mapsto f_t$ is in $\cS(\R_+^*;C_c^\infty(D(M_0,M)))$.

Now, let  $\theta :V'\to V$ be an ``exponential map'' which is a diffeomorphism of a (relatively compact) neighborhood $V'$ of the $0$ section $M_0$ in $\cN$ onto a tubular neighborhood $V$ of $M_0$ in $M$. We obtain a diffeomorphism $\Theta:W'\to W$ where $W'=\{(x,U,t)\in \cN\times \R;\ (x,tU)\in V'\}$ and $W=V\times \R^*\cup \cN\times \{0\}$.

Since $D(M_0,M)=M\times \R^*\cup W$, it follows that $C_c^\infty (D(M_0,M))=C_c^\infty (M\times \R^*)+C_c^\infty (W)$, and since both $M\times \R^*$ and $W$ are invariant by $\alpha$, it follows that $\cS_\alpha(D(M_0,M))$ is the sum of $\cS(\R^*;C_c^\infty (M))$ obtained as $\varphi_\alpha(\cS(\R_+^*;C_c^\infty(M\times \R^*)))$ and $\{f\in C^\infty(W);\ f\circ \Theta\in \cS_\alpha(\cN\times \R)\}$ (where $f\circ \Theta$ is extended by $0$ outside $W'$).

We denote by $\gJ_0( M_0,M)$ the subspace $\cS(\R^*;C_c^\infty (M))$ of $\cS_\alpha(D(M_0,M))$. \label{ectoire}

For $f\in \cS_\alpha(D(M_0,M))$ and $t\in \R$, we denote by $f_t:z\mapsto f(z,t)$ (with $z\in M$ if $t\ne 0$ and $z\in \cN$ if $t=0$). 

Note that if $f\in \cS_\alpha(D(M_0,M))$ and $\psi \in C^\infty (M)$ vanishes in the neighborhood of $M_0$, then the map $(z,t)\mapsto \psi (z)f(z,t)$ for $t\ne 0$ extends to a smooth map on $M\times \R$ vanishing at infinite order on $M\times \{0\}$.

As a direct consequence of Proposition \ref{J(G)1} we have:

\begin{proposition}\label{margot}
Let  $\theta :V'\to V$ be an exponential diffeomorphism as above and $\chi\in C_c^\infty(V)$ equal to $1$ near $M_0\subset M$. Let $f\in \cS_\alpha(D(M_0,M))$ \begin{enumerate}
\item For all $g\in C_c^{\infty}(\cN)$ the function $F:M_0\times \R\to \C$ defined by $F(x,0)=\hat f_0(x,0) g(x,0)$ and $F(x,t)=t^{-p}\int_{\cN_x} g(x,U)f_t\circ \theta (x,U)\,dU$ for $t\ne 0$  is smooth on $M_0\times \R$.\label{margota}
\item The following are equivalent:\label{gique}
\begin{enumerate}\renewcommand\theenumii{\roman{enumii}}
\renewcommand\labelenumii{\rm ({\theenumii})}
\item For all $g\in C_c^{\infty}(M)$ the function $t\mapsto \int_M g(x)f (x,t)\,dx$ vanishes as well as all its derivatives at $0$.
\item  For all $g\in C_c^{\infty}(\cN)$ the function $(x,t)\mapsto \int_{\cN_x} g(x,U)f_t\circ \theta (x,U)\,dU$ vanishes as well as all its derivatives on $M_0\times \{0\}\subset M_0\times \R$.
\item The function $(x,\xi,t)\mapsto \widehat {(\chi f_t)\circ \theta }(x,\xi)$ vanishes as well as all its derivatives on $M_0\times \{0\}$ sitting in $\cN^*\times \R$ as zero section.
\end{enumerate}

\end{enumerate}
In particular, it follows that conditions {\rm (ii)} and {\rm(iii)} do not depend on the choice of $\theta$. \\
We denote by $\gJ(M_0,M)\subset \cS_\alpha(D(M_0,M))$ the set of functions satisfying the above equivalent conditions.
\begin{proof}
For every $f\in \cS_\alpha(D(M_0,M))$, the function $(1-\chi )f_t$  has rapid decay when $t\to 0$. Replacing $f$ by $\chi f$, we may assume that $f\in \cS_\alpha(D(M_0,W))$. The result follows from Proposition \ref{J(G)1} since  $\cS_\alpha(D(M_0,W)){\buildrel {\theta^*}\over \longrightarrow}\cS_\alpha(D(M_0,W'))$ is an isomorphism.
\end{proof}
\end{proposition}

\subsection{A family of semi-norms}\label{seminormes}

We define a family $N_{k,\ell,j,m}$ of semi-norms on smooth functions on $\R^n\times \R^p\times \R$, for $k\in \N^n$, $\ell\in \N^p$, $j\in \N$ and $m\in \Z$. Put $$N_{k,\ell,j,m}(f)=\sup_{(x,\xi,t)\in \R^n\times \R^p\times \R}(\|\xi\|^2+t^2)^{m/2}\left |\frac{\partial ^{|k|+|\ell|+j}f}{\partial x^k\partial \xi ^\ell \partial t^j}(x,\xi,t)\right|.$$

We now use the notation introduced in previous subsection. Assume first that $M_0$ is compact. Fix a finite open cover $({\cal O}_i)_{i\in I}$ of $M_0$ by subsets diffeomorphic to $\R^n$ over which the normal bundle is trivialized. Using a partition of the identity adapted to $({\cal O}_i)_{i\in I}$, we obtain a family of semi-norms $N_{k,\ell,j,m}^i$ on $C^\infty(\cN^* \times \R)$ and thus, via the Fourier transform and the map $\Theta:W'\to W$ defined above, a family $\widetilde N_{k,\ell,j,m}^i$ of semi-norms on the space $C_c^\infty(W)$ defined by $$\widetilde N_{k,\ell,j,m}^i(f)=N_{k,\ell,j,m}^i(\widehat {f\circ \Theta})$$ (where $f\circ \Theta$ is extended by $0$ outside $W'$).

Finally, if $M_0$ is not compact, we define similarly a family of semi-norms using a locally finite cover $({\cal O}_i)_{i\in I}$.

 \emph{The action of $\R_+^*$}\label{actionR}

The action of $\R_+^*$ on $D(M_0,M)$ leads to an action by automorphisms $u\mapsto \alpha_u$ on $\cS_\alpha(D(M_0,M))$ and its subspaces $\gJ(M_0,M)$ (see prop. \ref{margot}) and $\gJ_0(M_0,M)$ ($=\cS(\R^*;C_c^\infty (M))$). Note that $\widehat {\alpha_u(f)\circ \Theta}(x,\xi,t)=u^p \widehat {f\circ \Theta}(x,u\xi,ut)$ and therefore, the semi-norms $\widetilde N_{k,\ell,j,m}^i$ are multiplied by a suitable power of $u$ by this action.

\section{The groupoids}

\subsection{Main notation} 

Let us recall  some standard constructions and notation on groupoids. 

\paragraph{Densities.}
If $M$ is a smooth manifold, $E$ is a real vector bundle of dimension $p$ on $M$ and $q\in \R$, we denote by $\Omega ^q(E)=|\Lambda ^pE^*|^q$ the $q$-density bundle on $E$: for $x\in M$, $\Omega^q_xE$ is the set of maps $\psi:\Lambda^pE_x\setminus\{0\}\to \R$ such that  $\psi(\lambda X)=\vert \lambda \vert^q \psi(X)$. A $q$-density is a section of $\Omega^q(E)$ and the product of a $q$-density with a $q'$-density leads to a $(q+q')$-density. There is a natural isomorphism $\Omega^q(E\oplus E')\to \Omega ^q(E)\otimes \Omega^q(E')$.

The positivity of a density makes sense, and thus the bundle of densities is an oriented real line bundle which is therefore trivial(izable). It is however sometimes important to keep track of the natural normalizations they give rise to.

We put $\Omega^p(M)=\Omega^p(TM)$. 

The main use of $1$-densities  is that their  integral over $M$ makes sense.  The natural action of diffeomorphism takes into account the Radon-Nykodym derivative and therefore, there is a unique linear form $\int_M : C_c^\infty(M;\Omega^1(M)) \rightarrow \R$ which agrees in local coordinates with the Lebesgue integral. In this way one associates to a submersion $p:M\to M_1$ and a vector bundle $E$ on $M_1$ a natural map $p_!:C_c^{\infty}(M;p^*E\otimes \Omega^{1} (\ker dp))\to C_c^{\infty}(M_1;E)$ obtained by integrating one densities along the fibers of $p$.

\paragraph{Source, range, algebroid.}
When $\cG \rightrightarrows \cG^{(0)}$ is a Lie groupoid with source $s$ and range $r$, we denote $\cG_x:= s^{-1}(x)$ and $\cG^x:=r^{-1}(x)$ for any $x\in \cG^{(0)}$. The set of composable elements is $\cG^{(2)}=\{(\gamma,\gamma');\ s(\gamma)=r(\gamma')\}$; the product $(\gamma,\gamma')\mapsto  \gamma \gamma' $ is a smooth submersion $p:\cG^{(2)}\to \cG$.

The \emph{$s$-vertical tangent bundle} is the tangent space to the $s$-fibers, that is $T_s \cG := \ker ds= \underset{x\in \cG^{(0)}}{\cup} T\cG_x$. The $r$-vertical tangent bundle $T_r \cG := \ker dr$ is defined similarly. Recall that the restriction of $T_s \cG$ to the set $\cG^{(0)}$ of units identifies with the total space $\gA \cG$ of the Lie algebroid  of $\cG$, which can also be thought of as the normal bundle to the inclusion $\cG^{(0)}\to \cG$. As usual, we will denote by $T^*_s \cG$ and $\gA^* \cG$ the corresponding dual bundles. 

\paragraph{The $*$-algebra of a  groupoid.}
Let $\cG$ be a Lie groupoid. As explained in \cite{ConnesSurvey}, the natural $*$-algebra is obtained using half densities, namely $C_c^{\infty}(\cG;\Omega^{1/2}(\ker ds\oplus \ker dr))$ endowed with the following operations:

\begin{description}
\item[Involution.] The map $\kappa:\gamma\mapsto \gamma^{-1}$ exchanges $r$ and $s$ and therefore it acts naturally on $\Omega^{1/2}(\ker ds\oplus \ker dr)$; also, the bundle $\Omega^{1/2}(\ker ds\oplus \ker dr)$ has a real structure, \ie there is a natural complex conjugation $\omega \mapsto \overline\omega $ of this bundle. The adjoint of $f\in C_c^{\infty}(\cG;\Omega^{1/2}(\ker ds\oplus \ker dr))$ is defined by $f^*(\gamma)=\kappa_*\Big(\overline{f(\gamma^{-1})}\Big)$.

\item [Product.] If $f,g\in C_c^{\infty}(\cG;\Omega^{1/2}(\ker ds\oplus \ker dr))$, then the restriction of $f\otimes g$ to $\cG^2$ is a section of the bundle $\Omega ^{1/2}(p^*T_r \cG\oplus \ker dp\oplus \ker dp\oplus p^*T_s \cG)=p^*\Omega ^{1/2}(\ker dr\oplus \ker ds)\otimes \Omega^1\ker dp$, and by integration along the fibers of $p$ we obtain $f\ast g\in C_c^{\infty}(\cG;\Omega^{1/2}(\ker ds\oplus \ker dr))$.
\end{description}

From now on, we just write $C_c^{\infty}(\cG;\Omega^{1/2})$ instead of $C_c^{\infty}(\cG;\Omega^{1/2}(\ker ds\oplus \ker dr))$. The $C^*$-algebra of the groupoid $\cG$ is a completion of the $*$-algebra $C_c^{\infty}(\cG;\Omega^{1/2}))$.
\begin{description}
\item[The reduced $C^*$-algebra] is obtained as the completion of $C_c^{\infty}(\cG;\Omega^{1/2})$ by the family of representations $(\lambda_x)_{x\in \cG^{(0)}}$, where $\lambda_x$ is the representation by left convolution on $L^2(\cG_x)$  (which is the completion of $C_c^\infty (\cG_x;\Omega^{1/2}\cG_x)$).

\item[The full $C^*$-algebra] is obtained as the completion of $C_c^{\infty}(\cG;\Omega^{1/2})$ by the family of all continuous representations (\cf \cite{Renault, KhoshSkand}).
\end{description}

\paragraph{Pseudodifferential operators on Lie groupoids}

Let $\cG$ be a Lie groupoid and let  $\theta :V'\to V$ be an ``exponential map'' which is a diffeomorphism of a (relatively compact) neighborhood $V'$ of the $0$ section $\cG^{(0)}$ in $\gA \cG$ (considered as the normal bundle to the inclusion $\cG^{(0)}\subset \cG$) onto a tubular neighborhood $V$ of $\cG^{(0)}$ in $\cG$. We assume that $r(\theta(x,U))=x$ for $x\in \cG^{(0)}$ and $U\in \gA_x \cG$.

Le $m\in \Z$. A classical pseudo-differential operator of order $m$ is a multiplier $P=P_0+K$ of the $*$-algebra $C_c^{\infty}(\cG;\Omega^{1/2})$ where $K\in C_c^{\infty}(\cG;\Omega^{1/2})$ and left multiplication by $P_0$ is given by an expression $$P_0\ast f(\gamma)=\int _{\xi \in \gA^*_{r(\gamma)} \cG}\Big(\int_{\gamma_1\in \cG^{r(\gamma)}} e^{i\langle \theta^{-1}(\gamma_1)|\xi\rangle} \varphi(r(\gamma),\xi)\chi(\gamma_1) f(\gamma_1^{-1}\gamma) \,d\gamma_1\Big)d\xi$$ where $\chi \in C_c^{\infty}(V)$ is a bump function satisfying $\chi=1$ on $\cG^{(0)}$ and $\varphi\sim \sum_{k=0}^{+\infty}a_{m-k}$ is a polyhomogeneous symbol ($a_j(x,\xi)$ is homogeneous of order $j$ in $\xi$).

We may then write $P_0(\gamma_1)=\int _{\xi \in \gA^*_{r(\gamma_1)} \cG} e^{i\langle \theta^{-1}(\gamma_1)|\xi\rangle} \varphi(r(\gamma_1),\xi)\chi(\gamma_1) \,d\xi$ meaning \eg that, as a multiplier, $P_0$ is the limit when $R\to \infty$ of $P_0^R\in C_c^{\infty}(\cG;\Omega^{1/2})$ where $$P_0^R(\gamma)=\int _{\xi \in \gA^*_{r(\gamma)} \cG;\ \|\xi\|\le R} e^{i\langle \theta^{-1}(\gamma )|\xi\rangle} \varphi(r(\gamma ),\xi)\chi(\gamma ) \,d\xi.$$ 

Recall that if the order $m$ of $P$ is strictly negative, then $P$ extends to an element of $C^*(\cG)$ and if $m=0$, then $P$ extends to a multiplier of $C^*(\cG)$. The map $P\mapsto a_0$ is well defined and extends to an onto morphism $\sigma_0:\Psi^*(\cG)\to C(S^*\gA \cG)$ with kernel $C^*(\cG)$ where $\Psi^*(\cG)$ is the closure of the algebra of pseudodifferential operators of order $\le 0$ in the multiplier algebra of $C^*(\cG)$ and $C(S^*\gA \cG)$ is the commutative $C^*$-algebra of continuous functions on the space $S^*\gA \cG$ of half lines in $\gA ^*\cG$. In other words, we have an exact sequence of $C^*$-algebras $$0\to C^*(\cG)\longrightarrow \Psi^*(\cG) {\buildrel{\sigma_0}\over\longrightarrow } C(S^*\gA \cG)\to0.$$

\subsection{The adiabatic groupoid}

\label{adiagrou}Let us start with a smooth groupoid $G \rightrightarrows G^{(0)}$ with source map $s$ and range map $r$ and denote by $\gA G$ its Lie algebroid: it is the normal bundle of the inclusion $G^{(0)}\to G$ as unit space.  

Its adiabatic groupoid  is the deformation to the normal cone $G_{ad}=D(G^{(0)},G)$ as in section \ref{DNC}.

As a set,  $G_{ad}=G\times \R^*\cup \gA G\times \{0\}$, and its set of objects $G_{ad}^{(0)}$ is $G^{(0)}\times \R$. The inclusions $q\mapsto (q,\lambda)$ are groupoid morphisms from $G$ to $G_{ad}$ for $\lambda \ne 0$, and from $\gA G$ to $G_{ad}$ for $\lambda=0$.

We fix an everywhere positive smooth density $\omega$ on $G$. We then get a smooth density $\omega_{ad} $ on $G_{ad}$: for $t\ne 0$, $\omega_{ad,t} =|t|^{-p}\omega$. We will now on fix the density $\omega$ (and therefore $\omega_{ad}$) and consider all elements of the groupoid algebra $C_c^\infty (G_{ad})$ as functions.

\subsection{The action of $\R_+^*$ and the \BIG}

The action of the group $\R^*_+$ on the deformation to a normal cone that we already used, is compatible with the groupoid structure of $G_{ad}=D(G^{(0)},G)$: 
$$\begin{array}{ccl} G_{ad}\times \R^*_+ & \rightarrow & G_{ad} \\   (\gamma,t;\lambda) & \mapsto & (\gamma,\lambda t) \mbox{ when } t\not=0 \\ (x,U,0,\lambda) & \mapsto & (x,\frac{1}{\lambda} U,0) \end{array}$$

The \emph{\BIG} is the (smooth) groupoid obtained as a crossed product of this action (see \eg \cite{LG}). Namely $\boinG:= G_{ad} \rtimes \R^*_+ \rightrightarrows G^{(0)}\times \R_+$ with structural morphisms 
$$\begin{array}{ll} \mbox{source and target :} & s_{bi}(\gamma,t;\lambda)=(s(\gamma), t) \mbox{ and } r_{bi}(\gamma,t;\lambda)=(r(\gamma),t\lambda ) \mbox{ for } t\not= 0 \\  & s_{bi}(x,U,0;\lambda)=r_{bi}(x,U,0;\lambda)=(x,0) \\ \mbox{product :} & (\gamma,\lambda' t;\lambda)\cdot(\gamma ' , t; \lambda')=(\gamma \gamma',t;\lambda \lambda')  \mbox{ for } t\not= 0  \\ &  (x,U,0;\lambda)\cdot (x,U',0;\lambda')=(x,\lambda' U+U',0;\lambda \lambda')\ . \end{array} $$

\begin{remarks}
\begin{enumerate}
\item At the level of Lie algebroids, these constructions are very simple and natural. Let us denote $(\gA G, \sharp,[\ ,\ ]_{\gA G})$ the Lie algebroid of $G$, with its corresponding anchor and bracket. 

The Lie algebroid of $G_{ad}$ is $(\gA G \times \R, \sharp_{ad},[\ ,\ ]_{ad})$ where $\sharp_{ad}:\gA G \times \R\to TG^{(0)}\times T\R$ is defined by $\sharp_{ad}(x,U,t)=(\sharp(x,tU),(t,0))$ for $(x,U,t)$ in $\gA G \times \R$,  and $[\ ,\ ]_{ad}$ is the Lie bracket which satisfies $[X,Y]_{ad}(x,t)=t[X,Y]_{\gA G}(x)$ where $X,\ Y$ are smooth (local) sections of $\gA G$.

 The Lie algebroid of $\boinG$ is $(\gA G \times T\R,\sharp_{bi},[\ ,\ ]_{bi})$ where $\sharp_{bi}(x,U,t,\lambda)=(\sharp(x,tU),(t,t\lambda))$ and $[\ ,\ ]_{bi}$ is the Lie bracket  induced by $[(X,\tau),(Y,\sigma)]_{bi}(x,t)=t([X,Y]_{\gA G},[\tau,\sigma])$ where $X,\ Y$ are smooth (local) sections of $\gA G$ and $\tau,\ \sigma$ are smooth local vector fields on $\R$.

\item This construction immediately extends to the case where $G$ is only assumed to be longitudinally smooth (\ie a continuous family groupoid in the sense of A. Paterson \cite{Paterson}).

\item Also, if $(M,\cF)$ is a singular foliation  in the sense of  \cite{AndrSk1} generated by vector fields $(X_i)_{1\le i\le n}$,  the \emph{adiabatic foliation} was constructed in \cite{AndrSk3}. We may construct the \emph{\BIF} on $M\times \R$ to be the foliation generated by the vector fields $t(X_i\otimes 1)$ and $t (1\otimes \partial /\partial t)$ (where $t$ is the $\R$ coordinate in $M\times \R$).
\end{enumerate}
\end{remarks}

\section{Schwartz algebra and module}

\subsection{The Schwartz algebra of Carrillo-Rouse and the ideal $\cJ(G)$}

Let $G$ be a Lie groupoid and $G_{ad}=D(G^{(0)},G)$ the corresponding adiabatic groupoid

We will use the Schwartz space described above for a general deformation to the normal cone; this leads to a slight modification of the Schwartz algebra of P. Carrillo-Rouse (\cf \cite{CR}).

\begin{defnot}\label{helene}
\begin{enumerate}
\item \emph{The ideal $\cJ_0(G)$ of functions with rapid decay at $0$} 

This is the space $\gJ_0(G^{(0)},G)=\cS(\R^*;C_c^\infty (G))$ defined above (section  \ref{renee} and  \ref{ectoire}).

It consists of smooth half densities $f$ with compact support on the groupoid $G\times \R_+$ such that, for every $k\in \N$, the function $(\gamma,t)\mapsto t^{-k}f(\gamma,t)$ extends smoothly on $G\times \R_+$. 

\item \emph{The (modified) Schwartz algebra $\cS_c(G_{ad})$ of Carrillo-Rouse}

This is the algebra $\cS_\alpha(D(G^{(0)},G))$ defined above (section \ref{ectoire}).

Its elements are sums $f+g$ where $f\in \cJ_0(G)$ and $g\in C^\infty(W;\Omega^{1/2})$ is such that $\widetilde N_{k,\ell,j,m}^i(g)<+\infty$ for all $i,k,\ell,j,m$ with $m\ge 0$.

Note that there is a canonical groupoid morphism $G_{ad}\to G$ (the image of $\gA G\times \{0\}$ is $G^{(0)}\subset G$) and, by definition, the image under this morphism of the support of $f\in \cS_c(G_{ad})$ is compact in $G$. On the other hand, unlike the original definition in \cite{CR}, we had to drop here the conical support requirement. 

\item \emph{The ideal $\cJ(G)$}

This is the space $\gJ(G^{(0)},G)$ defined in prop. \ref{margot}.

Its elements are sums $f+g$ where $f\in \cJ_0(G)$ and $g\in C^\infty(W;\Omega^{1/2})$ is such that $\widetilde N_{k,\ell,j,m}^i(g)<+\infty$ for all $i,k,\ell,j,m$ with any $m\in \Z$.
\end{enumerate}
\end{defnot}

We now check that $\cS_c(G_{ad})$ is indeed an algebra and $\cJ(G)$ an ideal. We begin with a remark:

\begin{remark} An element $f\in \cS_c(G_{ad})$ is a family $(f_t)_{t\in \R}$, where $f_t\in C^\infty_c(G)$ for $t\ne0$ and $f_0\in \cS(\gA G)$. Note that, since $G\times \R^*$ is dense in $G_{ad}$,  $f$ is determined by $(f_t)_{t\ne0}$. Roughly speaking, the definition implies that the support of $f_t$  concentrates around $G^{(0)}$ when $t$ goes to $0$.
\end{remark}

\begin{proposition}
The space $\cS_c(G_{ad})$ is a $*$-algebra: for $f,g\in \cS_c(G_{ad})$, the families $(f_t^*)_{t\in \R}$ and $(f_t\ast g_t)_{t\in \R}$ belong to $\cS_c(G_{ad})$. Moreover $\cJ_0(G)$ is a $*$-ideal of $\cS_c(G_{ad})$.
\begin{proof}
The function $(\gamma_1,\gamma_2)\mapsto f(\gamma_1)g(\gamma_2)$ defined on $G_{ad}^{(2)}=D(G^{(0)},G^{(2)})$ is an element of $\cS_\alpha(D(G^{(0)},G^{(2)}))$. Since the composition $G^{(2)}\to G$ is a submersion, equivariant with respect to $\alpha$ (since $\alpha$ is an action of $\R_+^*$ by groupoid automorphisms), we find that integration along the fibers yields a continuous map $\cS_\alpha(D(G^{(0)},G^{(2)}))\to \cS_\alpha(D(G^{(0)},G))$.

By continuity of the product of $C_c^\infty(G)$, we find that if $f\in \cJ_0(G)$ or $g\in \cJ_0(G)$, then $f\ast g\in \cJ_0(G)$.

The assertions about the $*$-operations are obvious.
\end{proof}
\end{proposition}

\begin{lemma}\label{convolution1} Let $f\in \cS_c(G_{ad})$. \begin{enumerate}
\item For every $g\in C_c^{\infty}(G)$ the function $F:G\times \R\to \C$ defined by $F(\gamma,t)=f_t\ast g(\gamma)$ for $t\ne 0$ and $F(\gamma,0)=\hat f_0(r(\gamma),0)  g(\gamma)$ is smooth.
\item \label{convolution2}We have $f\in \cJ(G)$ if and only if, for any $g\in C_c^\infty(G)$, the family $(f_t\ast g)_{t\in \R^*}$ is an element of $\cJ_0(G)$.
\end{enumerate}
\begin{proof}
Let  $\theta :V'\to V$ be an ``exponential map'' which is a diffeomorphism of a (relatively compact) neighborhood $V'$ of the $0$ section $G^{(0)}$ in $\gA G$ onto a tubular neighborhood $V$ of $G^{(0)}$ in $G$. We assume that $r(\theta(x,U))=x$ for $x\in G^{(0)}$ and $U\in \gA_x G$. 
\\
Let $\chi \in C_c^\infty (G)$ with support in $V$, such that  $\chi(\gamma)=1$ for $\gamma$ near $G^{(0)}$.  Then $((1-\chi )f_t)_{t\in \R}\in  \cJ_0(G)$, whence $((1-\chi )f_t\ast g)_{t\in \R}\in \cJ_0(G)$ for all $g\in C_c^\infty(G)$. 
\\
Furthermore, $$(\chi f_t\ast g)(\gamma)=\int_{\gA_{r(\gamma)} G} f_t\circ \theta^{-1}(r(\gamma),U)h(\gamma,U)\,dU,$$ where $h(\gamma,U)=\chi(\theta (r(\gamma),U))g(\theta (r(\gamma),U)^{-1}\gamma)\delta (r(\gamma),U)$ - here $\delta$ is a suitable Radon-Nykodym derivative.
\\
The Lemma follows now from Prop. \ref{margot}.
\end{proof}
\end{lemma}

\begin{proposition}
The space $\cJ(G)$ is a $*$-ideal of the algebra $\cS_c(G_{ad})$.
\begin{proof}
Follows immediately from Lemma \ref{convolution1}.\ref{convolution2}).
\end{proof}
\end{proposition}

\subsection{The smooth module $\cE^\infty$} 

The following rather technical Lemma uses the action $\alpha$ of $\R_+^*$ and the semi-norms  defined in section \ref{seminormes} and used in the Definitions \ref{helene}.

\begin{lemma}\label{teknik}
Let $f,g\in \cJ(G)$. We assume that their support is small enough (in $G$) so that, for every $t,u\in \R^*$, the function $f_t\ast g_u$ has support in $V$. The function $u\mapsto (f_t\ast g_{tu})_{t\in \R}$ is a smooth map from $\R_+^*$ to $\cJ(G)$ with rapid decay as $u\to 0$ and as $u\to\infty$. More precisely, for every $i\in I$, $k\in \N^n$, $\ell\in \N^p$, $j\in \N$, $m\in \Z$ and $q\in \Z$, the function $u\mapsto u^q\widetilde N_{k,\ell,j,m}^i(f\ast\alpha _u(g))$ is bounded on $\R_+^*$.

\begin{proof} 
We may perform the construction of the adiabatic groupoid starting from the adiabatic groupoid $G_{ad}$! Since $g=(g_t)$ is an element of $\cJ(G)$, the function $(t,u)\mapsto \chi (t)\chi(u) g_{tu}$ is an element of  $\cJ(G_{ad})$ where $\chi\in C_c^\infty(\R)$ is equal to $1$ near $0$. By Lemma \ref{convolution1} applied to the groupoid $G_{ad}$, it follows that $(\chi(t)\chi(u)f_{t}\ast g_{tu})_{u\in \R_+^*}$  leads to an element of $\cJ_0(G_{ad})$ and thus has rapid decay when $u\to 0$, uniformly in $t$. In other words, the functions $u\mapsto u^{-q}\widetilde N_{k,\ell,j,m}^i(f\ast \alpha _u(g))$ are bounded for  every $i\in I$, $k\in \N^n$, $\ell\in \N^p$, $j\in \N$, $m\in \N$ and $q\in \N$.

In the same way, $u\mapsto \alpha _u(f)\ast g$ has rapid decay when $u\to 0$, whence, using the compatibility of the semi-norms with the action of $\R_+^*$,  $u\mapsto f\ast \alpha ^{-1}_u(g)$ has rapid decay when $u\to 0$.  From this, we deduce that the functions $u\mapsto u^{q}\widetilde N_{k,\ell,j,m}^i(f\ast \alpha _u(g))$ are bounded for  every $i\in I$, $k\in \N^n$, $\ell\in \N^p$, $j\in \N$, $m\in \N$ and $q\in \N$, whence $q\in \Z$.

In other words, $u\mapsto f\ast \alpha _u(g)$ has rapid decay from $\R_+^*$ to $\cS_c(G_{ad})$. To see that it also has rapid decay as a function from $\R_+^*$ to $\cJ(G)$, it is enough to check that, for $h\in C_c^\infty(G)$, the map $u\mapsto f\ast \alpha _u(g)\ast h$ has rapid decay from $\R_+^*$ to $\cJ_0(G)$. But by Lemma \ref{convolution1}.\ref{convolution2}), $(t,u)\mapsto f_t\ast g_{tu}\ast h$ is an element of $\cS(\R^2;C_c^\infty(G))$ and vanishes as well as all its derivatives when $t=0$ or $u=0$.
\end{proof}\end{lemma}

\begin{theorem}\label{pseudosintegrales}
For $f\in \cJ(G)$ and $m\in \N$, the operator $\displaystyle\int_0^{+\infty}\! t^m f_t\, \displaystyle\frac {dt}{t}$ is an order $-m$ pseudodifferential operator of the groupoid $G$ \ie an element of $\cP_{-m}(G)$; its principal symbol $\sigma$ is given by $\sigma(x,\xi)=\displaystyle\int_0^{+\infty} t^m\hat f(x,t\xi,0)\displaystyle\frac {dt}{t}\cdot$ More precisely, there is a (classical) pseudodifferential operator $P$ with principal symbol $\sigma$ such that for every $g\in C_c^\infty (G)$, we have $P\ast g=\displaystyle\int_0^{+\infty}\! t^m f_t\ast g\, \displaystyle\frac {dt}{t}$ and $g\ast P=\displaystyle\int_0^{+\infty}\! t^mg\ast f_t\, \displaystyle\frac {dt}{t}\cdot$

\begin{proof}
Of course, if $f\in \cJ_0(G)$ then $\int_0^{+\infty} t^m f_t\,\frac{dt}t \in C_c^\infty (G)$. In particular, this gives the meaning of  $\displaystyle\int_0^{+\infty}\! t^m f_t\ast g\, \displaystyle\frac {dt}{t}$ and $\displaystyle\int_0^{+\infty}\! t^m g\ast f_t\, \displaystyle\frac {dt}{t}$ (thanks to Lemma \ref{convolution1}).

We thus need only to treat the case of $f\in C_c^\infty(W;\Omega^{1/2})$ satisfying the above condition. We may trivialize the half densities using a positive half density $\omega $ on $G$; we then may write $f_t=t^{-p}h_t \omega$ where $p$ is the dimension of the fibers of $G$ and $h$ is the restriction to $G\times \R_+^*$ of a smooth function with compact support in $G_{ad}$.  Such a function can be written as $h(\gamma,t)=\chi(\gamma)\chi'(t)\varphi\Big(\frac{\theta^{-1} (\gamma)}t,t\Big)$, where $\chi $ and $\chi'$ are bump-functions $\chi\in C_c^\infty(G)$ with support contained in $V$ which is equal to $1$ in a neighborhood of $G^{(0)}$, $\chi'\in C_c^\infty(\R_+)$ equal to $1$ in a neighborhood of $0$ and $\varphi\in C_c^\infty(\gA G\times \R_+)$ such that $\hat \varphi$ vanishes as well as all its derivatives at points of the form $(x,0,0)$ with $x\in G^{(0)}$.

Writing $\varphi(x,X,t)=(2\pi)^{-p}\int e^{i\langle X|\xi\rangle}\hat \varphi(x,\xi ,t)\,d\xi $, we find 
\begin{eqnarray*}
f_t(\gamma)&=&(2\pi t)^{-p}\chi(\gamma)\chi'(t)\omega \int e^{i\langle \frac{\theta^{-1} (\gamma)}t|\xi\rangle}\hat \varphi(x,\xi,t)\,d\xi\\
&=&(2\pi )^{-p}\chi(\gamma)\chi'(t)\omega \int e^{i\langle \theta^{-1} (\gamma)|\xi\rangle}\hat \varphi(x,t\xi,t)\,d\xi.
\end{eqnarray*}
 Therefore, we have an equality (as multipliers of $C_c^\infty (G;\Omega^{1/2})$), $$\int_0^{+\infty} t^mf_t(\gamma)\,\frac {dt}t=(2\pi)^{-p}\chi(\gamma)\omega\int e^{i\langle \theta^{-1} (\gamma)|\xi\rangle}a(x,\xi)\,d\xi,$$ where 
 \begin{eqnarray*}
a(x,\xi)&=&\int_0^{+\infty} t^m\chi'(t)\hat \varphi(x,t\xi ,t)\,\frac {dt}t\cdot\end{eqnarray*}

Taking derivatives of $a$ in $x$ gives the same type of expression; taking derivatives in $\xi$ increases $m$. 

An expression $c(x,\xi)=\int_0^{+\infty}  (t\|\xi\|)^m\chi'(t)b(x,t\xi ,t)\,\frac {dt}t$ with $b$ having rapid decay at infinity and points $(x,0,0)$ is bounded.
In other words, $a$ is a symbol of order $-m$ and type $(1,0)$

Now, given an expansion $\hat \varphi(x,\zeta,u)\sim \sum_{k=0}^\infty b_k(x,\zeta)u^k$, (for $u$ small) we find an expansion (for $\xi$ large), $a(x,\xi )\sim \sum_{k=0}^\infty a_{k+m}(x,\xi)$ where $a_{k+m}(x,\xi)=\int_0^\infty b_k(x,t\xi)t^{k+m}\,\frac {dt}t$ is homogeneous of degree $-k-m$ in $\xi$.
\end{proof}\end{theorem}

Of course, the same computation works for $m\in \Z$ and even $m\in \C$ (\cf \cite{chief} for pseudodifferential operators on groupoids with complex order).

\medskip Note that  Theorem \ref{pseudosintegrales}  shows that every element $P\in \cP_0(G)$ can be written as an integral $P_f=\displaystyle\int_0^{+\infty}\! f_t\, \displaystyle\frac {dt}{t}$ with $f\in \cJ(G)$ (using the standard Borel's Theorem-like techniques to construct $f$ such that $P-P_f\in C_c^\infty(G)$; of course if $P\in C_c^\infty(G)$, then $P=P_f$ where $f_t=\chi(t)P$ with an obvious choice of $\chi\in \cS(\R_+^*)$).

\begin{lemma}\label{ActionADroite}Let $f=(f_t)_{t\in \R}\in \cJ(G)$ and $P\in \cP_0(G)$.  There is a (unique) element $h=(h_t)_{t\in \R}$ of $\cJ(G)$ such that $h_t=f_t\ast P$ for $t\in \R^*$. Moreover $\widehat {h_0}=\widehat {f_0}\sigma_0(P)$. 

\begin{proof}
Uniqueness follows from density of $G\times \R^*$ in $G_{ad}$.

The family $(f_t\ast P)$ is smooth in $G\times \R^*$ and the image in $G$ of its support is compact.

Of course, if $f\in \cJ_0(G)$ then $(f_t\ast P)\in \cJ_0(G)$, we may thus assume that $f$  has support in a neighborhood of $G^{(0)}$ in $G$ as small as we wish.

In the same way, since $P$ is quasi-local, we may assume,  thanks to Lemma \ref{convolution1}, that the support of $P$ is also contained in a suitable neighborhood of $G^{(0)}$ in $G$.

Now, using Theorem \ref{pseudosintegrales}, we may write $P=\int g_u \frac{du}{u}$ with $g\in \cJ(G)$.

By locality, we may assume that, for all $t,u\in \R^*$, the support of $f_t\ast g_u$ is contained in $V$.

We find $h_t=f_t\ast P=\int_0^{+\infty}f_t\ast g_u\frac{du}{u}=\int_0^{+\infty}f_t\ast g_{tu}\frac{du}{u}$. The conclusion follows immediately from Lemma \ref{teknik}.

We found $h=\int_0^{+\infty} f\ast \alpha_u(g)\frac{du}{u}$. Evaluating at $t=0$, we find $h_0=\int_0^{+\infty}f_0\ast \alpha_u(g_{0})\frac{du}{u}$, and thus $\widehat{h_0}=\widehat {f_0}\int_0^{+\infty} \alpha_u(\widehat{g_{0}})\frac{du}{u}=\widehat {f_0}\sigma_0(P)$ (by Theorem \ref{pseudosintegrales}).
\end{proof}
\end{lemma}

Lemma \ref{ActionADroite} asserts that $\cJ(G)$ is endowed with a (right) $\cP_0(G)$-module structure; we will denote this smooth module by $\cE^{\infty}$.

\section{The Hilbert module $\cE$}

From now on, we will be interested in the restriction $G_{ad}^+$ of $G_{ad}$ to $G^{(0)}\times \R_+$. Let ${\rm ev_+}:C^*(G_{ad}) \rightarrow C^*(G_{ad}^+)$ be the morphism induced by the restriction map. Since $G^{(0)}\times \R_+$ is a closed saturated subspace of $G_{ad}$ and $G_{ad}^+$ is invariant under the action of $\R^*_+$, all the previous results obviously remain true when one replaces $\cS_c(G_{ad})$, $\cJ(G)$ and $\cJ_0(G)$ by their image under ${\rm ev_+}$. For the simplicity of notations we keep the same notation for $\cJ(G)$ and $\cJ_0(G)$: they are the image under ${\rm ev_+}$ of the previously defined  $\cJ(G)$ and $\cJ_0(G)$. 

\subsection{Completion of $\cE^\infty$}

Let $\Psi^*(G)$ denote the $C^*$-algebra of pseudodifferential operators, \ie the norm closure of $\cP_0(G)$ in the multiplier algebra of $C^*(G)$. Let also $\sigma_0:\Psi^*(G)\to C(S^*\gA G)$ be the principal symbol map. We have the exact sequence of $C^*$-algebras $$0\to C^*(G)\longrightarrow \Psi^*(G){\buildrel{\sigma_0}\over\longrightarrow } C(S^*\gA G)\to 0\eqno(2)$$
For $P\in \Psi^*(G)$, the function $\sigma_0(P)$ is thought of as a homogeneous function defined outside the zero section in $\gA G$.

The elements of $C^*(G_{ad}^+)$ are families $(f_t)_{t\in \R_+}$ with $f_t\in C^*(G)$ for $t\ne 0$ and $f_0\in C^*(\gA G)\simeq C_0(\gA^* G)$.

Put $J_0(G)=\{f\in C^*(G_{ad}^+);\  f_0 =0\}\simeq C_0(\R_+^*;G)$. We have an exact sequence $$0\to J_0(G)\longrightarrow C^*(G_{ad}^+){\buildrel{{\rm ev}_0}\over\longrightarrow } C_0(\gA^* G) \to 0.$$  Finally, put $J(G)\subset C^*(G_{ad}^+)=\{f\in C^*(G_{ad}^+);\ \forall x\in G^{(0)},\  \widehat {f_0}(x,0)=0\}$.

\begin{lemma}\label{essentialideal}
An element $f\in C^*(G_{ad}^+)$ is determined by the family $(f_t)_{t\in \R_+^*}$: the ideal $J_0(G)$ is essential  in $C^*(G_{ad}^+)$.
\begin{proof}
Let $x\in G^{(0)}$. Put $G_{ad,x}^+=\{(\gamma,t)\in G_{ad}^+;\ s(\gamma)=x\}=G_x\times \R_+^*\cup \gA_xG$. The map $(\gamma,t)\to t$ being a submersion, we obtain a continuous family $(H_{x,t})_{t\ge 0}$ of Hilbert spaces with $H_{x,0}=L^2(\gA_xG)$ and $H_{x,t}=L^2(G_x)$ for $t>0$ as a completion of smooth half densities with compact support on $G^+_{ad,x}$. 

Let $f\in C^\infty_c(G_{ad}^+;\Omega^{1/2})$ and $g\in C^\infty_c(G^+_{ad,x};\Omega^{1/2})$; we have $f* g\in C^\infty_c(G^+_{ad,x};\Omega^{1/2})$. It follows that $t\mapsto \|f_t* g_t\|$ is continuous; by density, this remains true for $f\in C^*(G^+_{ad})$ and any continuous section $g=(g_t)$ of $(H_{x,t})$. But $\|f_0\|=\sup\{\|f_0* g_0\|,\ x\in G^{(0)};\ g\in (H_{x,t});\ \sup_{t} \|g_t\|\le 1\}$. It follows that $\|f_0\|\le \sup_{t>0} \|f_t\|$.
\end{proof}
\end{lemma}

In the following, we consider $C^*(G^+_{ad})$ in the multiplier algebra $\cM(J_0(G))$ of $J_0(G)$. Since $J_0(G)=C_0(\R_+^*)\otimes C^*(G)$, the algebras $C^*(G)$, and $\Psi^*(G)$ sit also in $\cM(J_0(G))$.

From lemmas \ref{convolution1} and \ref{ActionADroite}, we immediately get.

\begin{proposition} \label{ActionADroite2}
 \begin{enumerate}
\item For $f\in C^*(G)$ and $g\in J(G)$ we have $f\ast g\in J_0(G)$. \label{convolution2}
\item For $P\in \Psi_0^*(G)$ and $f\in J(G)$ we have $f\ast P \in J(G)$ and $\widehat{(f\ast P)_0}=\widehat{f_0}\sigma_0(P)$. \label{ActionADroite2b}
\hfill $\square$\end{enumerate}
\end{proposition}

\begin{lemma}\label{strictConv}
For  $f\in \cJ(G)$, the integral $\displaystyle\int_0^{+\infty}\! f_t\, \displaystyle\frac {dt}{t}$ of Theorem \ref{pseudosintegrales} converges strictly (\ie in the topology of multipliers of $C^*(G)$).
\begin{proof} 
Let $g\in C_c^\infty(G)$. It follows from lemma \ref{convolution1} that $\int_0^{+\infty} f_t\ast g\,\frac{dt}t$ converges in norm. Taking adjoints, it follows that $\int_0^{+\infty} g\ast f_t\,\frac{dt}t$ converges also in norm. 

From Theorem \ref{pseudosintegrales}, it follows that $\int_0^{+\infty} f_t\ast g\,\frac{dt}t=P\ast g$ where $P$ is a pseudodifferential operator of order $0$ and therefore extends to a multiplier of $C^*(G)$. 

Now, if $f$ is a positive element in $\cJ(G)$, it follows that $\langle g|\Big(\displaystyle\int_s^{+\infty}\! f_t\, \displaystyle\frac {dt}{t}\Big)\ast g\rangle \le \langle g|P\ast g\rangle $, therefore the family $\Big (\displaystyle\int_s^{+\infty}\! f_t\, \displaystyle\frac {dt}{t}\Big )_{s>0}$ is bounded. It follows that  $\Big(\int_s^{+\infty} f_t\,\frac{dt}t\Big)\ast g$ converges to $P\ast g$ and $g\ast \Big(\int_s^{+\infty} f_t\,\frac{dt}t\Big)$ converges to $g\ast P$ for all $g\in C^*(G)$ (when $s\to 0$). Therefore $\displaystyle\int_0^{+\infty}\! f_t\, \displaystyle\frac {dt}{t}$ converges strictly.

By the polarization identity, we find that $\int_0^{+\infty} g_t^*\ast h_t\,\frac{dt}t$ converges in the multiplier algebra of $C^*(G)$ for   $g,h\in \cJ(G)$.

Now, one can find $g,h\in \cJ(G)$ such that $f_0=g_0^*\ast h_0$ (take for instance $\hat g_0(x,\xi)=|\hat f_0(x,\xi)|^2+\exp(-\|\xi\|^2-\|\xi\|^{-2})^{1/4}$ and $\hat h_0(x,\xi)=\hat f_0(x,\xi)\hat g_0(x,\xi)^{-1}$). It follows,  that there exists $f^1\in \cJ(G)$ with $f_t=tf^1_t+g_t^*\ast h_t$.  Now, since $f^1\in C^*(G_{ad}^+)$, it follows that $\|f_1^t\|$ is bounded. Therefore,  (using rapid decay at $\infty$) the integral $\int_0^{+\infty} f^1_t\,dt$ is norm convergent in $C^*(G)$.
\end{proof}
\end{lemma}

\begin{lemma}\label{fullmodule}
There exists $f\in \cJ(G)$ such that $\int _0^{+\infty} f_t^*\ast f_t\frac {dt}{t}$ is an invertible element of $\Psi^*(G)$ and $\Big(1-\int _0^{+\infty} f_t^*\ast f_t\frac {dt}{t}\Big)\in C_c^\infty(G)$.
\begin{proof}
Fix a smooth function  $\psi:\R_+\to \R_+$  with support in $]1,2[$ such that $\int_0^{+\infty}\psi^2(t)\frac {dt}{t}=1$. Let $g\in \cJ(G)$ be such that $\hat g_0(x,\xi)=\psi(\|\xi\|)$ and $g_t=0$ for $t\ge 1$. The positive pseudodifferential operator $P=\int _0^{+\infty} g_t^*\ast g_t\frac {dt}{t}$ has principal symbol equal to $1$ by Theorem \ref{pseudosintegrales}. By \cite{chief}, there exists $Q\in \cP_0(G)$ such that $1-Q^*PQ\in C_c^{\infty}(G)$. Using the exact sequence (2), it follows that there exists $b\in C^*(G)$ such that $b^*b+Q^*PQ\ge 1$. By density of $C_c^{\infty}(G)$ in $C^*(G)$, there exists $h\in C_c^{\infty}(G)$ such that $h^*\ast h+Q^*PQ$ is invertible. Taking $f_t=g_t\ast Q$ for $t\le 1$ and $f_t=\psi(t) h$ for $t \ge 1$, we find $\int _0^{+\infty} f_t^*\ast f_t\frac {dt}{t}=Q^*PQ+h^*\ast h$.
\end{proof}\end{lemma}

We may now construct the main object of this section.

\begin{theorem}\label{auboisdormant}
There is a Hilbert $\Psi^*(G)$-module $\cE$ containing $\cJ(G)$ as a dense subset, with the following operations:\begin{itemize}
\item  For $f,g\in \cJ(G)\subset \cE$, we have $\langle f|g\rangle =\int _0^{+\infty}f_t^*\ast g_t\frac{dt}t$. 
\item For $f\in \cJ(G)\subset \cE$ and $P\in \cP_0(G)\subset \Psi^*(G)$, we have $f\ast P\in  \cJ(G)\subset \cE$ and $(f\ast P)_t=f_t\ast P$ for $t\not=0$.
\end{itemize}
The module $\cE$ is a full $\Psi^*(G)$ module.
\begin{proof}
Clearly, for fixed $f$, the map $g\mapsto\langle f|g\rangle= \int _0^{+\infty}f_t^*\ast g_t\frac{dt}t$ is linear and $\cP_0(G)$-linear; also $\langle g|f\rangle=\langle f|g\rangle^*$. Furthermore $\langle f|f\rangle$ is the strict limit of elements of $C^*(G)_+$; therefore $\langle f|f\rangle\in \Psi^*(G)_+$. 

For $f\in \cJ(G)$, we then may put $\|f\|_\cE=\|\langle f|f\rangle\|_{C^*(G)}^{1/2}$. By the Cauchy-Schwarz inequality for $C^*$-modules, this defines a norm on $\cE^\infty$. Now, using again the Cauchy-Schwarz inequality, for $f,g\in \cJ(G)$ we have  $\|\langle f|g\rangle\|_{C^*(G)}\le \|f\|_\cE\|g\|_\cE$ and, for $P\in \cP_0(G)$, since $\langle f\ast P|f\ast P\rangle =P^*\langle f|f\rangle P\le \|\langle f|f\rangle\|_{C^*(G)}P^*\ast P$, we find $\|f\ast P\|_\cE\le \|f\|_\cE\|P\|_{\Psi^*(G)}$. It follows that the scalar product and the right action of $\cP_0(G)$ extend and endow $\cE$ with the desired Hilbert $\Psi ^*(G)$-module structure.

It follows from Lemma \ref{fullmodule} that $\cE$ is full - and in fact, that there exists $\zeta\in \cE$ ($\zeta=f\langle f|f\rangle^{-1/2}$) such that $\langle \zeta|\zeta\rangle=1$.
\end{proof}
\end{theorem}

\subsection{Computation of $\cK(\cE)$}

We now construct the desired natural isomorphism   $J(G)\rtimes \R_+^*\to \cK(\cE)$.

Note first, that if $f\in \cS_c(G^+_{ad})$ and $g\in \cJ(G)$, $f\ast g\in \cJ(G)$; furthermore, this left action is $\cP_0(G)$-linear and $\langle f\ast g|f\ast g\rangle=\int g_t^*\ast f_t^*\ast f_t\ast g_t\frac {dt}t\le \|f\|^2\langle g|g\rangle$ (where $\|f\|=\sup \|f_t\|$ is the norm of $f$ in $C^*(G^+_{ad})$ - this holds as well for the reduced and the full $C^*$-norm on $G$ and $G_{ad}$). Extending by continuity, we obtain a natural morphism $\pi_0:C^*(G^+_{ad})\to \cL(\cE)$.

The action of $\R_+^*$ on $G^+_{ad}$ gives rise to a unitary action on $\cE$ given by $(U_s(f))_t=f_{st}$ for $f\in \cJ(G)$ and $s,t\in \R_+^*$.

The couple $(\pi_0,U)$ is an equivariant representation of $(C^*(G^+_{ad}),\R_+^*)$, and therefore gives rise to a representation of $C^*(G^+_{ad})\rtimes \R_+^*\simeq C^*(\boinG)$, \ie a morphism $$\pi:C^*(G^+_{ad})\rtimes \R_+^*\to \cL(\cE) .$$

Put $\cE_0=\cE C^*(G)$. It is a closed submodule of $\cE$. 

Recall that if $J$ is a closed two sided ideal in  a $C^*$-algebra $A$ and $E$ is a Hilbert-$A$-module, then $EJ=\{x\in E;\ \langle x|x\rangle \in J\}$ is a closed submodule. The quotient $E/EJ$-is the Hilbert $A/J$-module $E\otimes _AA/J$. We clearly\footnote{Considering for instance the case $E=H_A$.} have a short exact sequence $$0\to \cK(EJ)\to \cK(E)\to \cK(E/EJ)\to0.$$

Note that if $J$ is an essential ideal of $A$, it follows that $EJ^\perp =\{0\}$. Indeed, for $x\in E$ non zero, there exists $b\in J$ such that  $\langle x|x\rangle b\ne 0$, whence $\langle x|xb\rangle =\langle x|x\rangle b\ne 0$. 

\begin{theorem}\label{bleue}
The morphism $\pi :J(G)\rtimes \R_+^*\to \cL(\cE)$ is an isomorphism from $J(G)\rtimes \R_+^*$ onto $\cK(\cE)$.
\begin{proof}
Consider the exact sequences:
$$\begin{array}{ccccccccc}
0&\to &J_0(G)\rtimes \R_+^*&\to &J(G)\rtimes \R_+^*&\to&(J(G)/J_0(G))\rtimes \R_+^*&\to &0\\
&&\pi \downarrow&&\pi \downarrow&&\overline\pi \downarrow&&\\
0&\to &\cK(\cE_0)&\to &\cK(\cE)&\to&\cK(\cE/\cE_0)&\to &0
\end{array}$$
We will show that:\begin{enumerate}
\item $\pi$ induces an isomorphism from $J_0(G)\rtimes \R_+^*$ onto $\cK(\cE_0)$;
\item $\pi(J(G)\rtimes \R_+^*)\supset \cK(\cE)$;
\item  The map $\overline \pi$ induced by $\pi$ gives rise to an isomorphism from $(J(G)/J_0(G))\rtimes \R_+^* $ onto $\cK(\cE/\cE_0)$.
\end{enumerate}
The Theorem then follows by diagram chasing.

We proceed with the proof of these facts:
\begin{enumerate}
\item According to Lemma \ref{convolution1}, the module $\cE_0$ is the closure of $\cJ_0(G)$. It is therefore canonically isomorphic to $C_0(\R_+^*)\otimes C^*(G)$ and, since $\R_+^*$ acts by translation in $C_0(\R_+^*)$ the statement follows.

\item Let $f,g\in \cJ(G)$; by Lemma \ref{teknik}, the map $s\mapsto f\ast \alpha_s(g^*)$ has rapid decay and thus defines an element $\int f\ast \alpha_s(g^*) \lambda _s ds/s$ in $J(G)\rtimes \R_+^*$. A direct computation then shows that $\vartheta _{f,g}=\pi\Big(\int f\alpha_s(g^*) \lambda _s ds/s\Big)$ where $\vartheta _{f,g}$ is the usual ``rank one'' operator $h\mapsto f\langle g|h\rangle$ on $\cE$. The result follows from density of $\cJ(G)$ in $\cE$.

\item The quotient $\Psi_*(G)/C^*(G)$ is isomorphic via the principal symbol map $\sigma$ to $C(S^* \gA G)$. Theorem \ref{pseudosintegrales} and Lemma \ref{ActionADroite} give the computations of $\sigma (\langle f|g\rangle)$ and $\widehat{(f\ast P)_0}$ for $f,g\in \cJ(G)$ and $P\in \cP_0(G)$. It follows that $\cE/\cE_0\simeq \cE\otimes_\sigma C(S^* \gA G)$ is the Hilbert $C(S^* \gA G)$ module $C(S^*\gA G)\otimes L^2(\R_+^*)$ obtained as completion of $C_c^{\infty}(\gA^*G\setminus G^{(0)})$ with respect to the $C(S^* \gA G)$ valued scalar product given by $\langle f|g\rangle(x,\xi)=\int_0^{+\infty}\overline{f(x,t\xi)}g(x,t\xi)\frac {dt}{t}$ and right action $(fh)(x,\xi)=f(x,\xi)h(x,\frac{\xi}{\|\xi\|})$ (for $f,g\in C_c^{\infty}(\gA^*G\setminus G^{(0)})$ and $h\in C^\infty(S^* \gA G)$). Moreover, left action of $J(G)/J_0(G)\simeq C_0(\gA^*G\setminus G^{(0)})$ is given by pointwise multiplication and the action of $\R_+^*$ is by scaling. The result follows. \qedhere
\end{enumerate}
\end{proof}
\end{theorem}

It is worth noting that thanks to Lemma \ref{essentialideal} and the amenability of $\R_+^*$, the ideal $J_0(G)\rtimes \R_+^*$ of $C^*(\boinG)=C^*(G^+_{ad})\rtimes \R_+^*$ is essential. In particular this remark gives an alternative proof of the injectivity of $\pi : J(G)\rtimes \R_+^*\to \cL(E)$.

\subsection{Pseudodifferential operators as convolution kernels}

Using Theorem \ref{bleue} together with Lemma \ref{fullmodule}, we can see $\Psi^*(G)$ as sitting as a corner in $C^*(G_{ga})$: Let $\zeta\in \cE$ satisfying $\langle \zeta|\zeta\rangle=1$ (as in Theorem \ref{auboisdormant}). Then $\vartheta_{\zeta,\zeta}\in \cK(E)\subset C^*(G_{ga})$  is a projection and $\Psi^*(G)\simeq \vartheta_{\zeta,\zeta}C^*(G_{ga})\vartheta_{\zeta,\zeta}$.

Actually, taking $f\in \cE^\infty$ such that $\langle f|f\rangle$ is invertible and $1-\langle f|f\rangle\in C_c^\infty(G)$ (Lemma \ref{fullmodule}), we can really see the elements of $\cP_0(G)$ as convolution operators on $G_{ga}$: the map $\vartheta_{T(f),f}\mapsto \langle f|T(f)\rangle$ is an isomorphism from a corner of the convolution algebra of smooth functions with Schwartz decay on $G_{ga}$ onto a subalgebra of $\cP_0(G)$, containing $\langle f|f\rangle \cP_0(G)\langle f|f\rangle$. Note that for every $P\in \cP_0(G)$, $P-\langle f|f\rangle\ast P\ast \langle f|f\rangle\in C_c^\infty (G)$. Finally $P=Q+R$ where $Q$ is (the image of) a smooth function on $G_{ga}$ and $R\in C_c^\infty(G)$.

\subsection{Stability of $\mathcal E$}

We now prove that the module $\cE$ is stable, \ie isomorphic to Kasparov's universal module $\cH_{\Psi^*(G)}$. 

We begin by recalling a few facts:

\begin{facts}
\begin{enumerate}
\item Recall first a few easy facts about Hilbert modules: \begin{enumerate}\renewcommand\theenumii{\arabic{enumii}}
\renewcommand\labelenumii{\rm {\theenumii}.}
\renewcommand\labelenumiii{\rm {(\theenumiii})}
\item A  Hilbert module $E$ is countably generated if and only if $\cK(E)$ is $\sigma$-unital (\ie has a countable approximate unit).
\item For a countably generated Hilbert module $E$ over a $C^*$-algebra $B$, the following are equivalent:\begin{enumerate}
\item the Hilbert $B$-module $E$ is stable  \ie isomorphic to $\ell^2(\N)\otimes E$;
\item the $C^*$-algebra $\cK(E)$ is stable \ie isomorphic to $\cK\otimes \cK(E)$;
\item there is a morphism $\psi:\cK\to \cL(E)$ such that $\psi(\cK)E=E$, or equivalently $\psi(\cK)\cK(E)=\cK(E)$.
\end{enumerate}
 Note that by Cohen's theorem (\cite{Cohen, HeRo} - see \eg \cite{Ped} for the case of $C^*$-algebras) there is no linear span or closure.
 
Indeed, implications (i) $\Rightarrow $ (ii) $\Rightarrow $ (iii) are straightforward. If (iii) is satisfied, then $E=\bigoplus _{k} \psi(e_{kk})E$ and the various $\psi(e_{kk})E$ are isomorphic \emph{via} $\psi(e_{jk})$ where $e_{jk}$ are matrix units for $\cK$, and (i) follows.

\item For  a countably generated Hilbert module $E$ over a $\sigma $-unital $C^*$-algebra $B$, the following are equivalent:\begin{enumerate}
\item $E$ is stable and full;
\item $E\simeq \cH_B$.
\end{enumerate}
Indeed (ii) $\Rightarrow $ (i) is obvious.

Conversely, let $E^*$ be the Hilbert $\cK(E)$-module $\cK(E,B)$; if $E$ is stable, then $\cK(E)$ is isomorphic as a $\cK(E)$-module to $\ell^2\otimes \cK(E)=\cH_{\cK(E)}$ and therefore $E'\oplus \cK(E)\simeq \cK(E)$ for every countably generated Hilbert $\cK(E)$-module $E'$;  if $E$ is full, then $B\simeq E^*\otimes_{\cK(E)} E$; finally, if (ii) is satisfied, then $$\cH_B\simeq \cH_B\oplus E\simeq ((\cH_B\otimes _BE^*)\oplus \cK(E))\otimes _{\cK(E)}E\simeq \cK(E)\otimes _{\cK(E)}E\simeq E.$$
\end{enumerate}

\item Let $\cG$ be a longitudinally smooth groupoid with compact space of objects $\cG^{(0)}$.

Every element of the algebroid of $\cG$, \ie a section of $\gA \cG$, defines a differential operator affiliated to $\cG$, \ie a multiplier of $\cG$ (once chosen a trivialization of the longitunal half densities on $\cG$).

Recall (\cf \cite{chief}) that the closure of an elliptic pseudodifferential operator $D$ on $\cG$ is a \emph{regular} unbounded multiplier of $C^*(\cG)$ and $(1+D^*\overline D)^{-1}$ is a strictly positive element of $C^*(\cG)$; moreover, if $D$ is formally self-adjoint, then $\overline D$ is self-adjoint. 

In particular, if $(X_1,\ldots ,X_m)$ are elements spanning the algebroid of $\cG$ as a $C^\infty(\cG^{(0)})$ module, then the closure of $\sum_i X_i^*X_i$ is such a self-adjoint elliptic operator whose spectrum is in $\R_+$ - \ie it is positive.

Furthermore, if $f$ is a smooth everywhere positive function on $\cG^{(0)}$, the operator $\sum_i X_i^*X_i+f=f^{1/2}\Big(\sum \tilde X_i^*\tilde X_i+1\Big)f^{1/2}$  is invertible where $\tilde X_i=f^{-1/2}X_i$ and its inverse $\Big(\sum_i X_i^*X_i+f\Big)^{-1}$ is a strictly positive element in $C^*(\cG)$.

\item Recall that a regular self-adjoint positive multiplier $D$ of a $C^*$-algebra $\cA$ with resolvent in $\cA$ defines a morphism $\pi_D:f\mapsto f(D)$ from $C_0(\R_+^*)$ to  $\cA$. Note that, for $t\in \R_+^*$, we have $\pi_{tD}=\pi_D\circ \lambda_t$ where $\lambda_t$ is the automorphism of $C_0(\R_+^*)$ induced by the regular representation. Since $t\mapsto \frac{t}{t^2+1}$ is a strictly positive element of $C_0(\R_+^*)$, it follows that $\pi_D(C_0(\R_+^*))\cA$ is the closure of $D(D^2+1)^{-1}\cA$.
\end{enumerate}
\end{facts}

\begin{proposition}\label{afer}
Let  $G$ be a Lie groupoid with compact $G^{(0)}$ and $G_{ad}$ its adiabatic groupoid; let $(Y_1,\ldots ,Y_m)$ span $\gA G$ as a module over $C^\infty(G^{(0)})$. Let $D_1=\Big(\sum_i Y_i^*Y_i+1\Big )^{1/2}$. There is a unique morphism $\psi:C_0(\R_+^*)\to C^*(G_{ad}^+)$ such that \begin{enumerate}\renewcommand\theenumi{\roman{enumi}}
\renewcommand\labelenumi{\rm ({\theenumi})}

\item $ev_1\circ \psi=\pi_{D_1}$ where $ev_1:C^*(G_{ad})\to C^*(G)$ is evaluation at $1$.
\item $\alpha_u\circ \psi=\psi\circ \lambda_u$ for all $u\in \R_+^*$.\label{propaferii}
\end{enumerate}
Moreover, \begin{enumerate}
\item $ev_0\circ \psi(f)=f\circ q^{1/2}\in C_0(\gA ^*G)$ where $q=\sum Y_i^2$ the $Y_i$'s being considered as functions on $\gA^* G$ (linear on each fiber).\label{propafera}
\item $\psi(C_0(\R_+^*))\subset J(G)$ and $\psi(C_0(\R_+^*))J(G)=J(G)$.\label{propaferb}
\end{enumerate}

\begin{proof} 
If $\psi $ satisfies (i) and (ii), then $(\psi (f))_t=f(tD_1)$ for $t\ne 0$. This shows immediately uniqueness of $\psi$.

Choose $u_0\ge 1$ and let $\cG$ be the restriction of $G_{ad}$ to $[0,u_0]$. Let $(X_1,\ldots, X_m)$ be the canonical extension of $(Y_1,\ldots ,Y_m)$ to $\gA \cG$. In particular for $u\ne 0$ we have $(X_i)_u=uY_i$. Put $D=\Big(\sum_i X_i^*X_i+h^2\Big )^{1/2}$, where $h$ is the function $(x,u)\mapsto u$ defined on $\cG^{(0)}=G^{(0)}\times [0,u_0]$. Note also $D_1=\Big(\sum_i Y_i^*Y_i+1\Big )^{1/2}$ its evaluation at $1$. For $u\in ]0,u_0]$, since $D_u=uD_1$, we have $ev_u\circ \pi_D=\pi_{D_1}\circ \lambda_{u}$.

Since the spectrum of $D_1$ is contained in $[1,+\infty[$, it follows that  for $f\in C_0(\R_+^*)$ and  $u\in \R_+^*$,  $\|\pi_{D_1}\circ \lambda_{u}(f)\|\le \sup \{|f(v)|;\ v\ge u\}$. In particular $\lim_{u\to+\infty} \pi_{D_1}\circ \lambda_{u}(f)=0$. It follows that there is an element $\psi_f\in C^*(G_{ad}^+)$ such that, for $u\ne 0$ we have $(\psi_f)_u=\pi_{D_1}\circ \lambda_{u}$ and $\psi_f$ restricted to $[0,u_0]$ is equal to $f(D)$. 

We thus get a homomorphism $\psi:f\mapsto \psi_f$ from $C_0(\R_+^*)$ to $C^*(G_{ad})$.

It satisfies (i), and since $\alpha_u\circ \psi(f)$ and $\psi\circ \lambda_u(f)$ coincide on $\R_+^*$, they are equal. Property (ii) follows.

Now, $D_0^2=q$ so we get property (\ref{propafera}).

By uniqueness, the morphism $\psi $ does not depend on the choice of $u_0$. 

It follows from property (\ref{propafera}) that, for all $f\in C_0(\R_+^*)$, since $f(0)=0$, $ev_0\circ \psi(f)$ vanishes on $M\subset \gA^*G$, whence $\psi(f)\in J(G)$.

Put $J'=\psi(C_0(\R_+^*))J(G)$.

As the restriction of $h$ to the $[\varepsilon,u_0]$ is invertible for every $u_0$ and $\varepsilon$, the restriction of $D(1+D^2)^{-1}$ to $[\varepsilon,u_0]$ is a strictly positive element; it follows that $J'$ contains the functions $\R_+^*\to C^*(G)$ with compact support. It therefore contains the ideal $C_0(\R_+^*)\otimes C^*(G)$ of $C^*(G^+_{ad})$. 

As the quotient $C^*(G^+_{ad})/C_0(\R_+^*)\otimes C^*(G)\simeq C_0(\gA G^*)$ is abelian, it follows that the right  ideal $J'$ is two sided. Finally $ev_0(J')=C_0(\gA G^*\setminus M) $  from  property (\ref{propafera}) since $\frac {q^{1/2}}{1+q}$ is a strictly positive element of $C_0(\gA G^*\setminus M) $.
\end{proof}
\end{proposition}

\begin{corollary}
The $C^*$-algebra $J(G)\rtimes \R_+^*\simeq \cK(\cE)$ is stable and therefore the Hilbert $\Psi^*(G)$-module is stable, \ie isomorphic to $\ell^2\otimes \Psi^*(G)$.
\begin{proof}
Indeed, by condition (\ref{propaferii}) in prop. \ref{afer}, $\psi$ induces a morphism $\hat\psi$ from $\cK=C_0(\R_+^*)\rtimes _\lambda \R_+^*$ to $J(G)\rtimes _\alpha \R_+^*$ and by (\ref{propaferb}), $\hat\psi(\cK)\big(J(G)\rtimes _\alpha \R_+^*)=J(G)\rtimes _\alpha \R_+^*$.
\end{proof}\end{corollary}

\begin{remark}
One can in fact show that $J(G)$ is isomorphic to a crossed product $\Psi^*(G)\rtimes_\beta \R$ in such a way that the $\R_+^*$ action on $J(G)$ is intertwined with the dual action $\hat \beta$. It follows that  $J(G)\rtimes \R_+^*=\Psi^*(G)\otimes \cK$ - we use the duality of $\R$ with $\R_+^*$ given by $(t,u)\mapsto u^{it}$.

The action $\beta$ is given by $\beta_t(P)=D_1^{it}PD_1^{-it}$. The operator $D_1^{it}$ is pseudodifferential of complex order $it$ (\cf \cite{chief}), therefore $\beta_t(P)$ is pseudodifferential of order $0$; it follows also from \cite{chief} that it has the same principal symbol as $P$.

We may embed $\Psi^*(G)$ in the multiplier algebra of $J(G)$ setting $P.(f_u)_{u\in \R_+^*}=(P*f_u)_{u\in \R_+^*}$ and $(f_u)_{u\in \R_+^*}.P=(f_u*P)_{u\in \R_+^*}$  thanks to prop. \ref{ActionADroite2}.\ref{ActionADroite2b}); furthermore, we have a one parameter group $(D^{it})_{t\in \R}$ in the multipliers of $J(G)$. As $D_u$ and $D_1$ are scalar multiples of each other, we find in this way a covariant representation of  the pair $(\Psi^*(G),\R)$. Associated to this covariant representation of $(\Psi^*(G),\R)$ is a morphism from $\Psi^*(G)\rtimes \R$ into the multiplier algebra of $J(G)$, but since  the image of $C^*(\R)\subset \Psi^*(G)\rtimes_\beta \R$ is contained in $J(G)$, we get a homomorphism $\varphi :\Psi^*(G)\rtimes_\beta \R\to J(G)$.

Note that the image of $\Psi^*(G)$ is translation invariant, \ie invariant by the extension $\overline\alpha_u$ of $\alpha_u$ to the multiplier algebra, and that $\overline\alpha_u(D^{it})=u^{it}D^{it}$. This shows that $\varphi$ is an equivariant morphism from $(\Psi^*(G)\rtimes_\beta \R,\hat\beta)$ to $(J(G),\alpha)$.

Now $\beta_t$ restricts to an action of $\R$ on $C^*(G)$, and according to prop. \ref{ActionADroite2}.\ref{convolution2}) it follows that $\varphi(C^*(G)\rtimes_\beta \R)$ is contained in the ideal $C_0(\R_+^*)\otimes C^*(G)$ of $J(G)$. As $D_1^{it}$ is a multiplier of $C^*(G)$, this crossed product is trivial. More precisely: if $t\mapsto w_t$ is a continuous homomorphism from a group $\Gamma$ to the multiplier algebra of a $C^*$-algebra $A$, the (full or reduced) crossed product $A\rtimes_{ad\, w}\Gamma$ is isomorphic to the (max or min) tensor product $A\otimes C^*(\Gamma)$ where $A$ and $\Gamma $ map to the multiplier algebra respectively by $a\mapsto a\otimes 1$ and $t\mapsto w_t\otimes \lambda_t$. It follows that $\varphi(C^*(G)\rtimes \R)=C_0(\R_+^*)\otimes C^*(G)$. Now, at the quotient level, the action $\beta$ becomes trivial on symbols: we thus obtain equality  $\varphi((\Psi^*(G)\rtimes_\beta \R)=J(G)$.
\end{remark}

\end{document}